\documentclass[12pt,a4paper]{article}
\usepackage[T1]{fontenc}
\usepackage{lmodern,amsmath,amssymb,amsthm,microtype}
\usepackage[margin=25mm]{geometry}
\usepackage[hidelinks]{hyperref}
\usepackage{enumitem,needspace}
\setlist{itemsep=3pt,topsep=5pt}
\newcommand{\N}{\mathbb N}
\newcommand{\B}{\mathcal B}
\newcommand{\C}{\mathcal C}
\newcommand{\res}{\mathord{\upharpoonright}}
\newcommand{\id}{\operatorname{id}}
\newcommand{\dom}{\operatorname{dom}}
\newtheorem{theorem}{Theorem}[section]
\newtheorem{lemma}[theorem]{Lemma}
\newtheorem{proposition}[theorem]{Proposition}
\newtheorem{corollary}[theorem]{Corollary}
\theoremstyle{remark}
\newtheorem{remark}[theorem]{Remark}
\newtheorem{problem}[theorem]{Problem}
\newtheorem*{countableproposition}{Proposition}
\newtheorem*{generalproblem}{Problem}
\hypersetup{pdftitle={Condensations onto G-delta Subsets of Cantor Cubes},pdfauthor={Wies\l aw \'Sliwa},pdfkeywords={condensation, Cantor cube, G-delta subset, Baire space, transfinite recursion}}
\title{Condensations onto $G_\delta$ Subsets of Cantor Cubes}
\author{Wies\l aw \'Sliwa\\[4pt]
\small Institute of Mathematics, Faculty of Exact and Technical Sciences\\
\small University of Rzesz\'ow, Pigonia 1, 35-310 Rzesz\'ow, Poland\\
\small \href{mailto:wsliwa@ur.edu.pl}{\texttt{wsliwa@ur.edu.pl}}}
\date{27 November 2026}
\begin{document}
\maketitle
\begin{abstract}
We prove in ZFC that every nonempty $G_\delta$-subset of the Cantor cube $2^{\omega_1}$ is a continuous bijective image of $\N^{\omega_1}$. A simultaneous countable-coordinate reduction yields a Polish parameter space with $G_\delta$-dense sections. A parameter-preserving homeomorphism places the excluded closed sets in zero slices associated with distinct coordinate blocks, and a transfinite construction preserves an exact product identity for the fibres. The proof does not use an ambient homeomorphism extension theorem. Consequences include the corresponding result for powers of noncompact zero-dimensional Polish spaces and a corollary for non-Archimedean function spaces. For $\kappa>\aleph_1$, the product of the main theorem with an elementary condensation treats nonempty $G_\delta$-subsets of $2^\kappa$ determined by at most $\aleph_1$ coordinates. We characterize continuous surjectivity onto $2^\kappa\setminus\{0\}$ by the existence of a closed discrete subset of $\N^\kappa$ of cardinality $\kappa$, and refine such surjections to have singleton fibres over points of nonempty support of size less than $\kappa$. Classical compactness results imply that even surjectivity fails at weakly compact cardinals. The condensation problem at $\aleph_2$ remains open here.
\end{abstract}


\section{Introduction}
A \emph{condensation} of a space $X$ onto a space $Y$ is a continuous bijection $X\to Y$, not necessarily a homeomorphism. The problem of finding coarser topologies with prescribed properties has a long history in general topology. In particular, Problem 1 in the Scottish Book, posed by Banach in 1935, asks which metrizable spaces admit condensations onto metrizable compact spaces; see Banakh and Plichko \cite{BanakhPlichko} for the historical background, including related questions of Alexandroff.

Several results illustrate the scope of this problem. Pytkeev \cite[Theorem 1, p.~831]{Pytkeev} showed that a non-$\sigma$-compact Borel subspace of a Polish space admits a condensation onto the Hilbert cube (see also \cite{Pavlov}). Banakh and Plichko \cite{BanakhPlichko} proved that every Banach space of density $\aleph_1$ admits a condensation onto the Hilbert cube. In the context of Cartesian products, Pavlov \cite{Pavlov} observed that $\mathbb Q^\kappa$ condenses onto $[0,1]^\kappa$ for every infinite cardinal $\kappa$, where $\mathbb Q$ has its usual topology. Brian \cite{Brian} and Brian and Miller \cite{BM} studied condensations of completely ultrametrizable spaces. More recently, Osipov and Pytkeev \cite{OsipovPytkeev} proved that every metrizable space of weight $\lambda$ with $\lambda^{\aleph_0}=\lambda$ admits a condensation onto a Banach space of weight $\lambda$.

Recent contributions include Osipov \cite{Osipov2024} on consistency results for condensations of absolute $\kappa$-Borel sets onto compacta, and Juh\'asz, van Mill and Soukup \cite{JMS2026} on restrictions on compact or separable condensation targets for locally compact spaces and groups. Our question concerns condensations onto prescribed subspaces of Cantor cubes rather than onto an unspecified space with a given property. Already in 1929 Sierpi\'nski \cite[Th\'eor\`eme, p.~18]{Sierpinski} characterized the subsets of the real line that are continuous bijective images of the space of irrational numbers: in modern terminology, they are precisely the nonempty Borel subsets in which every point is a condensation point. For Borel subspaces of arbitrary Polish spaces, the corresponding characterization was proved by Mauldin and Sarbadhikari \cite[Theorem 7, pp.~443--444]{MauldinSarbadhikari} and is also recorded by Komisarski, Michalewski and Milewski \cite[Lemma 12, p.~62]{KMM}: every nonempty locally uncountable Borel subspace of a Polish space is a continuous bijective image of Baire space. Here locally uncountable means that every nonempty relatively open subset is uncountable. Conversely, the Lusin--Souslin theorem \cite[Theorem 15.1, p.~89]{Kechris} and the fact that every nonempty open subset of Baire space has cardinality $2^{\aleph_0}$ show that these conditions are necessary. In particular, this characterizes all subspaces of $2^\omega$ onto which $\N^\omega$ condenses; among the nonempty $G_\delta$-subsets, precisely those without isolated points qualify.

For nonmetrizable Cantor cubes, related work concerns extending homeomorphisms between closed subspaces to the whole cube; see Shchepin and Valov \cite{SV,SVext}. On the source side, Karassev, Shchepin and Valov \cite{KSV} extend homeomorphisms between compact subsets of $\N^\tau$ to autohomeomorphisms of that product. These extension results concern prescribed homeomorphisms of subspaces, whereas our problem asks for a continuous bijection from one fixed product onto a prescribed, possibly noncompact, target. Our proof does not invoke those extension theorems. The elementary product condensation onto the full cube and the classical Polish-space characterization do not by themselves control the fibres over excluded closed sets in an uncountable cube.

We work in ZFC. Set
$\N=\{1,2,\ldots\}$, $\B=\N^{\N}$ and $\C=2^{\N}$, where
$2=\{0,1\}$ and $\N$ are discrete. All products have the ordinary
product topology. In particular, $\N^{\omega_1}$ is not the metrizable
space $\omega_1^\omega$. The symbol $0$ denotes the zero point of a
binary product. Empty products are singletons. For a partial function $p$
and a set $A$, the notation $p\res A$ means the restriction of $p$
to $\dom p\cap A$. The ordinal indices of
a tuple $x\in\B^\alpha$ start at $0$, while the digits of each
$x_\xi\in\B$ start at $1$.

There is a classical condensation $\B\to\C$ (see Lemma~\ref{lem:basic-condensation} for a short proof). Taking its powers and reindexing the products gives condensations $\N^\kappa\to 2^\kappa$ for every infinite cardinal $\kappa$. Although the preceding result on Borel sets already implies the countable open-set case, we include an elementary argument, since the construction will also be useful in the uncountable setting.

\Needspace{12\baselineskip}
\begin{countableproposition}
For every nonempty open subset $U$ of the Cantor space $2^\omega$, there exists a condensation
\[
 \N^\omega\longrightarrow U.
\]
\end{countableproposition}
\begin{proof}
Cover $U$ by a countable family of basic clopen cylinders contained in $U$. Successively subtract the preceding finitely many members of the cover and discard empty differences. This gives a partition of $U$ into a nonempty finite or countably infinite family $(U_i)_{i\in I}$ of nonempty clopen subsets of $2^\omega$. Each $U_i$ is homeomorphic to $2^\omega$ by the characterization of Cantor space. Partition one discrete coordinate of $\N^\omega$ into $|I|$ countably infinite sets. The resulting clopen pieces $(V_i)_{i\in I}$ are each homeomorphic to $\N^\omega$. For every $i\in I$, the classical condensation above gives a condensation $V_i\to U_i$. Combining these maps on the open partition of the domain gives a continuous bijection onto $U$.
\end{proof}

In the nonmetrizable case, nonempty $G_\delta$-subsets of Cantor cubes cannot have isolated points: every nonempty $G_\delta$-subset of an uncountable Cantor cube contains a copy of the entire cube through each of its points (fix the countably many coordinates used by basic neighborhoods inside the defining open sets). This leads naturally to the following problem.

\begin{generalproblem}
For which uncountable cardinals $\kappa$ does every nonempty $G_\delta$-subset $B$ of the Cantor cube $2^\kappa$ admit a condensation
\[
 \N^\kappa\longrightarrow B?
\]
\end{generalproblem}

\Needspace{7\baselineskip}
The main contribution is an affirmative ZFC answer for the first uncountable cardinal $\aleph_1$, with no additional assumptions on the prescribed nonempty $G_\delta$ target.

\begin{theorem}\label{thm:main-gdelta}
If $B$ is a nonempty $G_\delta$-subset of $2^{\omega_1}$, then there is
a condensation
\[\N^{\omega_1}\longrightarrow B.\]
\end{theorem}

Write $B=\bigcap_{n\in\N}U_n$ with each $U_n$ open. The proof has
three substantial parts. First, we choose countably many parameter
coordinates so that all nonempty sections of all the $U_n$ are
$G_\delta$-dense. Second, we use disjoint countable coordinate blocks
to transform the closed sets excluded by the $U_n$ into closed sets
$K_n$ lying in zero slices associated with distinct coordinate blocks. Third, we construct a continuous
map for which the fibre over a point outside $\bigcup_nK_n$ is a
singleton and the fibre over a point of that union is empty.

The whole fibre must remain a Cartesian product throughout the recursion. A collection of individually possible auxiliary tuples need not be jointly possible: for instance, the diagonal in $\{0,1\}^2$ projects onto both factors but omits two pairs. Controlling the projections alone would therefore not justify the simultaneous choices at a successor stage or the independent prefix conditions at a limit stage. Proposition~\ref{prop:countable-exceptions} establishes an exact fibre identity, with a factor $Z_\alpha$ for each unresolved exception and a singleton for each resolved one. At the final stage the former factor is empty, whereas the latter factors give a unique preimage.

\paragraph{Where the hypotheses enter.} Section~\ref{sec:reflection} uses countable-coordinate reflection and the Baire theorem on compact countably determined cylinders to obtain the Polish parameter $V$ (Proposition~\ref{prop:polish-parameter-reduction}); Lemma~\ref{lem:parameter-absorption} supplies the condensation onto $V$. Section~\ref{sec:normalization} uses its countable clopen base, the compact binary coordinate blocks, and fresh coordinates (the union of the blocks used before any stage below $\omega_1$ is countable) to normalize the excluded sets (Proposition~\ref{prop:simultaneous-normalization}). In Section~\ref{sec:recursion}, compactness of the omitted binary factors gives closed projections and their limit identities. Countability of $\alpha<\omega_1$ makes the parameter spaces at each proper stage zero-dimensional Polish, so Lemmas~\ref{lem:parameter-extension} and \ref{lem:exceptional-collapse} apply. It also makes $Z_\alpha$ nonempty, while $Z_{\omega_1}$ is empty by the impossibility of an injection from $\omega_1$ into $\N$. Finally, the supremum of the countably many stages at which the exceptions are resolved is still below $\omega_1$. The reflection results do not require the index cardinal to be regular, but the recursion uses these particular properties of $\omega_1$; Proposition~\ref{prop:paramobstruction} identifies an obstruction to extending the parameter lemma to general nonmetrizable bases. All auxiliary constructions needed for the main proof are given explicitly. Appendix~\ref{app:verification} gives a dependency and verification guide to the recursion.

The central contribution is Theorem~\ref{thm:main-gdelta}, supported by the reflection, normalization and exact-fibre construction. Lemma~\ref{lem:parameter-extension} supplies a reusable auxiliary tool. Sections~6 and~7 give consequences by product or homeomorphism arguments. In Section~8, Corollary~\ref{cor:smallcoordinates} is a product consequence of the main theorem, whereas Proposition~\ref{prop:extent} and Theorem~\ref{thm:singleton} concern surjectivity, not the full condensation problem. Corollary~\ref{cor:compact-obstruction} applies classical compactness results to show that the unrestricted assertion fails at weakly compact cardinals.

\section{Elementary Polish constructions}

\begin{lemma}\label{lem:characterizations}
Let $X$ be a nonempty zero-dimensional Polish space.
\begin{enumerate}
\item If $X$ is compact and has no isolated points, then $X\cong\C$.
\item If no point of $X$ has a compact neighborhood, then $X\cong\B$.
\end{enumerate}
\end{lemma}
\begin{proof}
For (1), fix a compatible metric on $X$. Recursively partition every
member of the preceding finite clopen partition into finitely many,
at least two, nonempty clopen sets of diameter at most $2^{-m}$ at
level $m$. Such a refinement exists by compactness and a clopen base;
any cell left unsplit can be split because there are no isolated
points. The partitions form a finitely branching tree, every node of
which has at least two children. Match the root with $\C$. If a node
has $r\ge2$ children and is matched with a binary cylinder, partition
that cylinder into $r$ cylinders using the suffixes
$\{1^j0:0\le j<r-1\}\cup\{1^{r-1}\}$. Every child cylinder has at least one
additional fixed digit. Along every branch, the nested clopen sets
in $X$ and the nested binary cylinders therefore have singleton
intersections. Matching these intersections gives a bijection
$\C\to X$. At each fixed tree level there are only finitely many
matched binary cylinders, so their word lengths have a finite
maximum. Agreement on sufficiently long prefixes forces the
images into a set of arbitrarily small diameter, so the map is
continuous. A continuous bijection from a compact space to a
Hausdorff space is a homeomorphism.

For (2), fix a compatible complete metric $d$ on $X$. Every nonempty
clopen $V\subseteq X$ is complete and noncompact. Consequently it is
not totally bounded. Choose $\varepsilon_V>0$ such that $V$ cannot be
covered by finitely many sets of diameter less than
$\varepsilon_V$. Given $\eta>0$, a countable clopen base supplies a
countable cover of $V$ by clopen sets contained in $V$ of diameter
less than $\min(\eta,\varepsilon_V)$. Subtract the preceding
finitely many cover members successively, and discard empty pieces.
The resulting clopen partition must be infinite and all its members
have diameter less than $\eta$.

Starting with $V_\varnothing=X$, construct nonempty clopen sets
$V_s$, $s\in\N^{<\omega}$, so that
$(V_{s^\frown n})_{n\in\N}$ partitions $V_s$ and its members have
diameter at most $2^{-(|s|+1)}$. For each $t\in\B$, completeness and
the diameter bounds give a unique point in
$\bigcap_{m\ge1}V_{t\res m}$. Every point of $X$ determines a unique
branch. The resulting bijection $f:\B\to X$ is continuous by the
diameter estimates, and $f([s])=V_s$ for every basic cylinder $[s]$.
It is open and hence is a homeomorphism.
\end{proof}

Every nonempty clopen subspace of $\B$ satisfies part (2). Indeed, a
basic cylinder is a closed copy of $\B$, and $\B$ is noncompact
because its projection onto the discrete space $\N$ is onto. A
compact neighborhood cannot contain such a closed noncompact copy.

\Needspace{11\baselineskip}
\begin{lemma}\label{lem:basic-condensation}
There is a condensation $c:\B\to\C$.
\end{lemma}
\begin{proof}
Let $S=\{0\}\cup\{1/n:n\in\N\}$ with its usual topology. Every
bijection $a:\N\to S$ is continuous. Hence $a^{\N}:\B\to S^{\N}$
is a continuous bijection. The space $S^{\N}$ is nonempty, compact,
metrizable and zero-dimensional. No point is isolated, since any
basic neighborhood leaves a coordinate free. Apply
Lemma~\ref{lem:characterizations}(1). This classical construction
also appears in the literature on condensations; see
\cite[Lemma 10, p.~5, of arXiv:1106.1629v1]{Miller}.
\end{proof}

\begin{lemma}\label{lem:parameter-absorption}
If $P$ is a nonempty zero-dimensional Polish space and $L$ is
countably infinite, then $\B$ condenses onto $P\times2^L$.
\end{lemma}
\begin{proof}
The space $P\times\B$ is nonempty, zero-dimensional and Polish. A
nonempty basic open set contains a subspace $\{p\}\times Q$ which
is closed in the whole product and homeomorphic to $\B$, with $Q$
a basic clopen cylinder. Thus the product has no compact
neighborhood at any point. Lemma~\ref{lem:characterizations}(2)
gives a homeomorphism $u:\B\to P\times\B$. Compose it with
$\id_P\times c_L$, where $c_L:\B\to2^L$ is obtained from
Lemma~\ref{lem:basic-condensation} by reindexing.
\end{proof}

\begin{lemma}\label{lem:finite-approximation}
Let $v:\B\to\C$ be continuous and let $m\in\N$. There is a
condensation $c_m:\B\to\C$ such that
\[
 c_m(t)\res m=v(t)\res m\quad\hbox{if }t(1)<m,
\]
where $\res m$ denotes the first $m$ digits.
\end{lemma}
\begin{proof}
For each binary word $s$ of length $m$, put
$H_s=\{t:t(1)<m,\ v(t)\res m=s\}$. Every nonempty $H_s$ is clopen
in $\B$ and is homeomorphic to $\B$. For each such $s$, choose a
nonempty proper clopen subset $D_s$ of the target cylinder $[s]$.
It is homeomorphic to $\C$, so condense $H_s$ onto $D_s$ using
Lemma~\ref{lem:basic-condensation}. The remaining source
$\{t:t(1)\ge m\}$ is a nonempty clopen Baire space. The remaining
target is a nonempty clopen Cantor space: there are finitely many
chosen $D_s$, each proper in its own cylinder. Condense the remaining
source onto this target. The finite clopen partition makes the
combined bijection continuous. This also covers $m=1$, when every
$H_s$ is empty.
\end{proof}

\begin{lemma}\label{lem:parameter-extension}
Let $P$ be zero-dimensional Polish, let $Z$ be closed in $P$, and let
$r:Z\times\B\to\C$ be continuous. There is a continuous extension
$q:P\times\B\to\C$ such that $q(p,\cdot):\B\to\C$ is a
condensation for every $p\in P\setminus Z$.
\end{lemma}
\begin{proof}
If $Z=P$, take $q=r$. If $Z=\varnothing$, use the fixed map $c$ in
the second variable. Otherwise choose a compatible metric on $P$
and nested countable clopen partitions $\mathcal P_m$, $m\ge0$,
with $\mathcal P_0=\{P\}$ and mesh at most $2^{-m}$ for $m\ge1$.
They are obtained by countable small-diameter clopen covers and
successive subtraction, refining the previous partition.

For $p\notin Z$, use its first partition cell $W$ disjoint from $Z$.
Such a cell exists because $Z$ is closed and the mesh tends to zero.
These first cells form a pairwise disjoint clopen cover of
$P\setminus Z$. If $W$ is at level $m(W)\ge1$, its parent meets
$Z$; choose $z_W$ in that intersection. Apply
Lemma~\ref{lem:finite-approximation} to $t\mapsto r(z_W,t)$ and
$m(W)$ to obtain a condensation $c_W$. Set $q(p,t)=r(p,t)$ on
$Z\times\B$ and $q(p,t)=c_W(t)$ for $p\in W$.

Only continuity on $Z\times\B$ needs checking. Let
$(p_j,t_j)\to(z,t)$ with $z\in Z$ and $p_j\notin Z$. Write $W_j$
for the first cell of $p_j$ and $m_j=m(W_j)$. For fixed $k$,
eventually $p_j$ belongs to the level-$k$ cell containing $z$.
Hence $m_j>k$, so $m_j\to\infty$. Eventually $m_j\ge2$, and
\[
 d(z_{W_j},z)\le 2^{-(m_j-1)}+d(p_j,z)\longrightarrow0.
\]
Also $t_j(1)=t(1)<m_j$ eventually. Thus $q(p_j,t_j)$ and
$r(z_{W_j},t_j)$ agree on their first $m_j$ digits, while continuity
of $r$ gives $r(z_{W_j},t_j)\to r(z,t)$. Terms with $p_j\in Z$
are dealt with directly by continuity of $r$. The domain is
metrizable, so this sequential check proves joint continuity.
\end{proof}

The second factor in Lemma~\ref{lem:parameter-extension} may be
replaced by any space homeomorphic to $\B$, by conjugating with a
fixed homeomorphism. We will use this observation explicitly.

For $\alpha\le\omega_1$, set $P_\alpha=\B^\alpha$ and
\begin{equation}\label{eq:Z-definition}
 Z_\alpha=\{x\in P_\alpha:x_\xi(1)>1\ (\xi<\alpha),\quad
 x_\xi(1)\ne x_\eta(1)\ (\xi<\eta<\alpha)\}.
\end{equation}
These sets are closed. The set $Z_0=P_0$ is a singleton. Each
$Z_\alpha$ for $\alpha<\omega_1$ is nonempty, by an injection
$\alpha\to\{2,3,\ldots\}$, whereas
\begin{equation}\label{eq:Z-empty}
 Z_{\omega_1}=\varnothing.
\end{equation}
The first digits record a partial injection of the ordinal index set into $\{2,3,\ldots\}$; the remaining digits leave room for the Baire-space constructions. The purpose is to retain nonempty exceptional fibres at every countable stage, but eliminate them at $\omega_1$. Restrictions take $Z_\delta$ into $Z_\alpha$ for $\alpha\le\delta$.
For every nonzero limit $\lambda\le\omega_1$,
\begin{equation}\label{eq:Z-limit}
 Z_\lambda=\{x\in P_\lambda:x\res\alpha\in Z_\alpha
                         \text{ for all }\alpha<\lambda\},
\end{equation}
because each defining requirement involves at most two block indices.
For countable $\alpha$, the inclusion
$Z_{\alpha+1}\subset Z_\alpha\times\B$ is proper: choosing the
new first digit to be $1$ puts a point outside $Z_{\alpha+1}$.

\begin{lemma}\label{lem:exceptional-collapse}
For every $\alpha<\omega_1$, there is a continuous map
$r_\alpha:Z_\alpha\times\B\to\C$ with
$r_\alpha^{-1}(0)=Z_{\alpha+1}$ and with a singleton fibre over
every point of $\C\setminus\{0\}$.
\end{lemma}
\begin{proof}
The nonempty space $W_\alpha=Z_\alpha\times\B$ is
zero-dimensional Polish and nowhere locally compact, by the last
Baire factor. By Lemma~\ref{lem:characterizations}(2), it and all
its nonempty clopen subspaces are homeomorphic to $\B$.

Its nonempty open subspace $W_\alpha\setminus Z_{\alpha+1}$ can be
partitioned into countably infinitely many nonempty sets $V_j$,
$j=0,1,\ldots$, clopen in $W_\alpha$. To see this, disjointify a
countable clopen cover of that complement. If only finitely many
nonempty pieces result, split one Baire piece into countably many
clopen cylinders. For $j\ge0$ put
\[
 D_j=\{z\in\C:z(i)=0\ (1\le i\le j),\quad z(j+1)=1\}.
\]
The $D_j$ partition $\C\setminus\{0\}$ and each is homeomorphic
to $\C$. Condense $V_j$ onto $D_j$, and put $r_\alpha=0$ on
$Z_{\alpha+1}$. The fibre claims and continuity on the $V_j$ are
immediate. At a point of $Z_{\alpha+1}$, the complement of the
clopen set $V_0\cup\cdots\cup V_{k-1}$ is a neighborhood on which
the first $k$ output digits are zero. This proves continuity there.
\end{proof}

\begin{remark}\label{rem:nonextendible}
No surjectivity of $Z_{\alpha+1}\to Z_\alpha$ is claimed or used.
For example, the element of $Z_\omega$ with $x_j(1)=j+2$,
$j<\omega$, exhausts all permitted first digits and has no
extension to $Z_{\omega+1}$. All remaining digits may be chosen
arbitrarily. The later proof uses nonemptiness at each countable
stage, not extendibility of each exceptional branch.
\end{remark}

\section{Countable-coordinate reflection and Polish parameters}\label{sec:reflection}

For a set $I$ and a finite partial function $s:I\to 2$, put
\[
 [s]_I=\{x\in2^I:s\subseteq x\}.
\]
We also use this notation when the domain of $s$ is countable. Such a
countably determined cylinder is a nonempty compact $G_\delta$-subset of
$2^I$. Every nonempty $G_\delta$-subset of $2^I$ contains one: for a point
in a countable intersection of open sets, choose a basic cylinder about
that point in each open set, and take the union of the corresponding
finite partial functions. A subset is called $G_\delta$-dense if it meets every nonempty
$G_\delta$-subset. Consequently, it is $G_\delta$-dense if and only
if it meets every countably determined cylinder.

\begin{lemma}\label{lem:simultaneous-reflection}
Let $I$ be uncountable, let $(U_n)_{n\in\mathbb N}$ be open subsets of
$2^I$, and let $C\subseteq I$ be countable. There is a countably infinite
set $A\subseteq I$ containing $C$ with the following properties.
\begin{enumerate}
\item If $p$ is a finite partial function with $[p]_I\subseteq U_n$ and
$\operatorname{dom}p\not\subseteq A$, and
$r=p\res(\operatorname{dom}p\cap A)$, then there is an
uncountable family $\mathcal R$ of finite partial functions $q$ such that
$[q]_I\subseteq U_n$, $q\res A=r$, and the nonempty finite sets
$\operatorname{dom}q\setminus A$, $q\in\mathcal R$, are pairwise disjoint.
\item For every $a\in2^A$ and every $n$, the section
\[
 (U_n)_a=\{z\in2^{I\setminus A}:(a,z)\in U_n\}
\]
is either empty or $G_\delta$-dense in $2^{I\setminus A}$.
\end{enumerate}
\end{lemma}

\begin{proof}
Let $\mathcal S_n$ be the family of finite partial functions $s:I\to2$
with $[s]_I\subseteq U_n$. Choose a countably infinite $A_0\supseteq C$.
Having chosen $A_m$, for every $n\in\mathbb N$ and every finite partial
function $r:A_m\to2$, choose a maximal subfamily
\[
 \mathcal M_{m,n,r}\subseteq
 \{s\in\mathcal S_n:r\subseteq s,
       \ \operatorname{dom}s\setminus\operatorname{dom}r\ne\varnothing\}
\]
whose nonempty sets
$\operatorname{dom}s\setminus\operatorname{dom}r$ are pairwise disjoint.
Such a family exists by Zorn's lemma; it may be empty. Define
\[
 A_{m+1}=A_m\cup
 \bigcup_{\substack{n\in\mathbb N,\ r\text{ finite on }A_m\\
                    |\mathcal M_{m,n,r}|\le\aleph_0}}
                 \ \bigcup_{s\in\mathcal M_{m,n,r}}\operatorname{dom}s,
 \qquad A=\bigcup_{m<\omega}A_m.
\]
There are only countably many pairs $(n,r)$ at each stage, so all these
sets are countable.

To prove (1), choose $m$ with $\operatorname{dom}r\subseteq A_m$. The
family $\mathcal M_{m,n,r}$ is uncountable. Otherwise every domain of one
of its members lies in $A_{m+1}\subseteq A$, whereas
$\operatorname{dom}p\setminus\operatorname{dom}r$ is a nonempty subset
of $I\setminus A$. Thus $p$ could be adjoined to this family, contrary to
maximality. Since its nonempty remainders are pairwise disjoint, at most
countably many meet $A$. Discard those members. The remaining uncountable
family has precisely the asserted properties.

For (2), suppose $(U_n)_a\ne\varnothing$ and choose
$p\in\mathcal S_n$ compatible with $a$. If
$\operatorname{dom}p\subseteq A$, the section is the whole complementary
cube. Otherwise put $r=p\res(\dom p\cap A)$ and use (1). Given any
countable $J\subseteq I\setminus A$ and any $v\in2^J$, one can choose
$q\in\mathcal R$ with
$(\operatorname{dom}q\setminus A)\cap J=\varnothing$. The partial
functions $a$, $v$ and $q$ are compatible. Any extension of their union
to $I$ belongs to $U_n$ and has the prescribed coordinates. Thus
$(U_n)_a$ meets every countably determined cylinder.
\end{proof}

\begin{lemma}\label{lem:enlargement-and-intersection}
\begin{enumerate}
\item Let $A\subseteq I$ be countable and suppose that every nonempty
section $U_a$ of an open $U\subseteq2^I$ is $G_\delta$-dense. If
$L\subseteq I\setminus A$ is countable, then every section
$U_{(a,b)}$, $b\in2^L$, is empty when $U_a$ is empty, and is
$G_\delta$-dense in $2^{I\setminus(A\cup L)}$ otherwise. In particular,
\[
 \pi_{A\cup L}(U)=\pi_A(U)\times2^L.
\]
\item A countable intersection of open $G_\delta$-dense subsets of a
Cantor cube is $G_\delta$-dense.
\end{enumerate}
\end{lemma}

\begin{proof}
For (1), fix a nonempty section $U_a$, a point $b\in2^L$, a countable
$J\subseteq I\setminus(A\cup L)$, and $v\in2^J$. Apply
$G_\delta$-density of $U_a$ to the countable assignment $b\cup v$ on
$L\cup J$. This proves the assertion about every enlarged section and
also the projection identity.

For (2), let $(W_n)$ be the open $G_\delta$-dense subsets and let $Q$ be
a countably determined cylinder. Each $Q\cap W_n$ is open and dense in
$Q$. Indeed, every nonempty relatively open subset of $Q$ contains a
countably determined cylinder, obtained by adjoining finitely many
conditions to those defining $Q$. It must meet $W_n$. The space $Q$ is
nonempty compact Hausdorff, so the Baire category theorem gives
$Q\cap\bigcap_n W_n\ne\varnothing$. This proves the claim.
\end{proof}

\begin{proposition}\label{prop:polish-parameter-reduction}
Let $I$ be uncountable and let
$B=\bigcap_{n\in\mathbb N}U_n\ne\varnothing$, where the $U_n$ are open
in $2^I$. There are a countably infinite $A\subseteq I$ and a nonempty
zero-dimensional Polish space $V\subseteq2^A$ without isolated points
such that
\[
 V=\pi_A(B)=\bigcap_{n\in\mathbb N}\pi_A(U_n),
\]
and, for every $a\in V$, all sections $(U_n)_a$ and the section $B_a$
are $G_\delta$-dense in $2^{I\setminus A}$. The set $V$ is itself a
$G_\delta$-subset of $2^A$. Moreover, $A$ can be written as
$A=A_0\sqcup L$, with $L$ countably infinite, so that
$V=V_0\times2^L$ for a nonempty zero-dimensional Polish
$V_0\subseteq2^{A_0}$.
\end{proposition}

\begin{proof}
Apply Lemma~\ref{lem:simultaneous-reflection} and write $A_0$ for the
resulting countable coordinate set. Put
\[
 V_0=\bigcap_n\pi_{A_0}(U_n).
\]
Product projections are open, so $V_0$ is $G_\delta$ in $2^{A_0}$.
Certainly $\pi_{A_0}(B)\subseteq V_0$, and this set is nonempty. If
$a\in V_0$, each $(U_n)_a$ is nonempty and $G_\delta$-dense. By
Lemma~\ref{lem:enlargement-and-intersection}(2), their intersection $B_a$
is $G_\delta$-dense and in particular nonempty. Thus
$V_0=\pi_{A_0}(B)$.

Choose a countably infinite $L\subseteq I\setminus A_0$ and put
$A=A_0\cup L$. Applying
Lemma~\ref{lem:enlargement-and-intersection}(1) to every $U_n$, and then
part (2) to the enlarged sections, gives
\[
 \pi_A(B)=\bigcap_n\pi_A(U_n)=V_0\times2^L=:V.
\]
It also gives all the asserted section properties. The space $V$ is a
nonempty $G_\delta$-subspace of a metrizable Cantor cube, so it is
zero-dimensional and Polish. The factor $2^L$ ensures that it has no
isolated points. More explicitly, every nonempty relatively open subset
of $V$ contains a product $O\times W$ with $O\ne\varnothing$ and $W$
a nonempty open subset of $2^L$, and hence is uncountable.
\end{proof}

\begin{lemma}\label{lem:fresh-coordinate-witness}
Let $I$ be uncountable, let $A\subseteq I$ be countable, let $U\subseteq
2^I$ be open, and let $a\in2^A$ have a $G_\delta$-dense section $U_a$.
For every countable $J\subseteq I\setminus A$ there exist an open
neighborhood $W$ of $a$ in $2^A$, a nonempty finite
$F\subseteq I\setminus(A\cup J)$, and $s\in2^F$ such that
\[
 \{x\in2^I:x\res A\in W,
                  \ x\res F=s\}\subseteq U.
\]
\end{lemma}

\begin{proof}
Apply Lemma~\ref{lem:simultaneous-reflection} to the single open set $U$
with the prescribed initial set $A\cup J$, and denote its resulting
countable set by $D$. The compact set
\[
 T=\{t\in2^D:t\res A=a\}
\]
is nonempty. For each $t\in T$, the $G_\delta$-density of $U_a$ gives a
point of $U$ extending $t$, since $D\setminus A$ is countable. Choose a
finite partial function $p_t$ with $[p_t]_I\subseteq U$ compatible with
that point, and set $r_t=p_t\res(\dom p_t\cap D)$. The cylinders
$[r_t]_D$ cover $T$. Choose a finite subcover, with the associated pairs
denoted by $(p_i,r_i)$, $1\le i\le m$.

For each $i$ with $\operatorname{dom}p_i\subseteq D$, take $q_i=p_i$.
For every other $i$, property (1) of
Lemma~\ref{lem:simultaneous-reflection}, applied to $D$, supplies an
uncountable family of conditions in $U$ with root $r_i$ on $D$ and
pairwise disjoint nonempty remainders outside $D$. Successively choose
$q_i$ from these families so that their remainders outside $D$ are also
disjoint from all remainders chosen earlier. Only finitely many
coordinates have to be avoided at each such choice. The roots $r_i$ may be incompatible; only the remainders outside $D$ are combined. Consequently
\[
 s=\bigcup_{i=1}^m
       q_i\res(\operatorname{dom}q_i\setminus D)
\]
is a consistent finite partial function whose domain $F$ misses $D$.
If $F$ is empty, adjoin one coordinate from $I\setminus D$ with an
arbitrary binary value. This preserves all the conditions below and
makes $F$ nonempty.

The clopen set $C=\bigcup_{i=1}^m[r_i]_D$ contains $T$. Since
$2^D\setminus C$ is compact, its projection onto $2^A$ is closed and
does not contain $a$. Let
\[
 W=2^A\setminus\pi_A(2^D\setminus C).
\]
Then $W$ is an open neighborhood of $a$ and
$\pi_A^{-1}(W)\subseteq C$ within $2^D$. If $x$ satisfies the two
conditions in the statement, its restriction to $D$ extends some $r_i$,
and its restriction outside $D$ extends the remainder of $q_i$.
Therefore $q_i\subseteq x$, so $x\in[q_i]_I\subseteq U$.
\end{proof}

\section{Simultaneous normalization of the excluded sets}\label{sec:normalization}

\begin{proposition}\label{prop:simultaneous-normalization}
Let $|I|=\aleph_1$, let $A\subseteq I$ be countable, and let
$V\subseteq2^A$ be a nonempty zero-dimensional second-countable space.
Let $(U_n)_{n\in\mathbb N}$ be open subsets of $2^I$ such that $(U_n)_a$
is $G_\delta$-dense for every $a\in V$ and every $n$. Set
\[
 X=V\times2^{I\setminus A},\qquad F_n=X\setminus U_n,
 \qquad \mathcal C=2^{\mathbb N}.
\]
There is a homeomorphism, preserving the parameter in $V$,
\[
 h:X\longrightarrow
 Y=V\times\mathcal C^{\omega_1}
        \times(\mathcal C^{\omega_1})^{\mathbb N}
\]
such that every $K_n=h(F_n)$ is closed in $Y$ and
\[
 K_n\subseteq
 \{(a,y,(z_j)_{j\in\mathbb N})\in Y:z_n=0\}.
\]
\end{proposition}

\begin{proof}
Let $T=I\setminus A$ and let
$\mathcal T=\mathbb N\times\omega_1\times\mathbb N$. Fix a bijection
$e:\omega_1\to\mathcal T$. We recursively choose pairwise disjoint
countably infinite sets $J_t\subseteq T$, $t\in\mathcal T$, in this
order. If $t=e(\xi)=(n,\alpha,m)$, the set
\[
 S_\xi=\bigcup_{\eta<\xi}J_{e(\eta)}
\]
is countable, because $\xi<\omega_1$. The indices $\alpha$ and $m$ provide the desired output indexing; the witness condition depends on $n$ and the already used coordinates, not on an additional property of $\alpha$ or $m$. For each $a\in V$, apply
Lemma~\ref{lem:fresh-coordinate-witness} to $U_n$ and the forbidden set
$S_\xi$. Obtain an open $W_a\subseteq2^A$, a nonempty finite
$F_a\subseteq T\setminus S_\xi$, and $s_a\in2^{F_a}$ such that
\[
 \{x:x\res A\in W_a,
                 \ x\res F_a=s_a\}\subseteq U_n.
\]
By second countability, choose a countable subcover of $V$ from the
sets $V\cap W_a$. The union of the supports $F_a$ associated with this
subcover is countable and misses $S_\xi$. Only now choose a countably
infinite $J_t\subseteq T\setminus S_\xi$ containing this union.

A zero-dimensional second-countable space has a countable clopen
partition subordinate to any countable open cover: first take a
countable clopen refinement $(C_j)$ that covers the space, and take the nonempty sets $C_j\setminus\bigcup_{i<j}C_i$. Each block has its own partition; no common refinement over all blocks is required. Let $\mathcal P_t$ be such a partition of
$V$ subordinate to the chosen cover. For each $P\in\mathcal P_t$,
choose an associated witness $s_a$ with $P\subseteq W_a$, and put
\[
 D_{t,P}=[s_a]_{J_t}\subseteq2^{J_t}.
\]
This is a nonempty proper clopen subset, since its finite support is
nonempty. It satisfies
\[
 \{x\in X:x\res A\in P,
                 \ x\res J_t\in D_{t,P}\}\subseteq U_n.
\]

Both $D_{t,P}$ and its complement are nonempty clopen subspaces of a
metrizable Cantor cube and are therefore homeomorphic to $\mathcal C$.
Choose a homeomorphism
\[
 \theta_{t,P}:2^{J_t}\longrightarrow2\times\mathcal C
\]
that maps $D_{t,P}$ onto $\{1\}\times\mathcal C$ and its complement
onto $\{0\}\times\mathcal C$. For $a\in P$, use this map on the block
$J_t$. The resulting map
\[
 \Theta_t:V\times2^{J_t}\longrightarrow V\times2\times\mathcal C,
 \quad (a,x)\longmapsto(a,\theta_{t,P}(x))\quad(a\in P)
\]
is a homeomorphism: its definition and its inverse are continuous on
the open partition $(P)_{P\in\mathcal P_t}$.

Put $R=T\setminus\bigcup_{t\in\mathcal T}J_t$. Applying all the block
maps, and leaving the coordinates in $R$ unchanged, gives a
parameter-preserving homeomorphism
\[
 X\longrightarrow V\times2^R\times\mathcal C^{\mathcal T}
                         \times2^{\mathcal T}.
\]
For clarity, continuity of this map follows coordinatewise from that
of the $\Theta_t$. The inverse reconstructs each block using
$\Theta_t^{-1}$ and is likewise continuous coordinatewise. Thus no
uniformity in the partitions $\mathcal P_t$ is required.

For a point $\varepsilon\in2^{\mathcal T}$, group its distinguished binary coordinates by the explicit reindexing
\[
 z_n(\alpha)(m)=\varepsilon(n,\alpha,m)
 \quad(n,m\in\N,\ \alpha<\omega_1).
\]
This is a homeomorphism $2^{\mathcal T}\to(\mathcal C^{\omega_1})^{\N}$, with inverse given by the same coordinate equality. In particular, vanishing of all distinguished coordinates with first index $n$ is equivalent to $z_n=0$. The remaining factor
$2^R\times\mathcal C^{\mathcal T}$ is homeomorphic to
$\mathcal C^{\omega_1}$ by a fixed coordinate bijection:
$|R|\le\aleph_1$ and $|\mathcal T|=\aleph_1$. These regroupings give
the homeomorphism $h$ in the statement.

If $x\in F_n$, then for every $t=(n,\alpha,m)$ its block $x\res J_t$
lies outside $D_{t,P}$, where $P$ is the member of $\mathcal P_t$
containing its parameter. Otherwise the displayed inclusion would
give $x\in U_n$. The distinguished binary coordinate associated with
every such $t$ is therefore zero. Hence $z_n=0$ on $h(F_n)$. Finally,
$F_n$ is closed in $X$, so $K_n=h(F_n)$ is closed in $Y$.
\end{proof}

\Needspace{24\baselineskip}
\section{A construction for countably many excluded sets}\label{sec:recursion}

\begin{proposition}\label{prop:countable-exceptions}
Let $V$ be a nonempty zero-dimensional Polish space, and suppose
there is a condensation $b:\B\to V$. Put
\[
 Y=V\times\C^{\omega_1}\times(\C^{\omega_1})^{\N}.
\]
For each $n\in\N$, let $K_n$ be closed in $Y$ and suppose that
\begin{equation}\label{eq:zero-slices}
 K_n\subseteq\{(v,a,(z_j)_{j\in\N})\in Y:z_n=0\}.
\end{equation}
Then there is a condensation
\[
 \N^{\omega_1}\longrightarrow Y\setminus\bigcup_{n\in\N}K_n.
\]
The sets $K_n$ are allowed to be empty or to intersect one another.
\end{proposition}

\begin{proof}
We give the preparation and the induction separately. Fix the
condensation $c:\B\to\C$ from
Lemma~\ref{lem:basic-condensation}.

\medskip\noindent\emph{Projections of the excluded sets.}
For $\alpha\le\omega_1$, define
\[
 Y_\alpha=V\times\C^\alpha\times(\C^\alpha)^{\N},
 \qquad K_{n,\alpha}=\pi_\alpha(K_n),
\]
where $\pi_\alpha:Y\to Y_\alpha$ preserves $V$ and restricts all
ordinal-indexed tuples to the indices below $\alpha$. For
$\alpha\le\delta\le\omega_1$, write
$\pi_\alpha^\delta:Y_\delta\to Y_\alpha$ for the corresponding
prefix projection.

Every $K_{n,\alpha}$ is closed. Indeed, the coordinates omitted by
$\pi_\alpha$ form a compact product, and projection along a compact
factor is closed. To recall the argument, if a point of the target
is outside the projection of a closed set, its entire compact fibre
lies in the open complement. Finitely many product neighborhoods
cover that fibre; intersecting their first-factor neighborhoods
gives a neighborhood still disjoint from the projection. In
particular, this step does not assume that $V$ or $K_n$ is compact.
The projections satisfy
\begin{equation}\label{eq:K-coherence}
 \pi_\alpha^\delta(K_{n,\delta})=K_{n,\alpha}.
\end{equation}
For every nonzero limit $\lambda\le\omega_1$,
\begin{equation}\label{eq:K-limit}
 K_{n,\lambda}=
 \bigcap_{\alpha<\lambda}
 (\pi_\alpha^\lambda)^{-1}(K_{n,\alpha}).
\end{equation}
Here is an explicit proof of the reverse inclusion. Fix a point
$y$ on the right, with $V$-coordinate $v$. The section $K_n(v)$
in the remaining compact coordinate product is closed, hence
compact. For each $\alpha<\lambda$, its points with the prescribed
$\alpha$-prefix form a nonempty closed subset, by membership in
$K_{n,\alpha}$. These subsets are nested. Compactness yields a
point with all the prescribed prefixes and therefore with prefix
$y$ of length $\lambda$. This proves \eqref{eq:K-limit}, also when
$\lambda=\omega_1$.

\medskip\noindent\emph{Auxiliary exceptional sets.}
Define closed sets $E_{n,\alpha}\subseteq Y_\alpha$ by recursion:
\begin{align}
 E_{n,0}&=Y_0=V,\label{eq:E-zero}\\
 E_{n,\alpha+1}
 &=\{y^+\in Y_{\alpha+1}:\pi_\alpha^{\alpha+1}(y^+)
          \in K_{n,\alpha},\quad z_n^+(\alpha)=0\},
       \label{eq:E-successor}\\
 E_{n,\lambda}
 &=\bigcap_{\alpha<\lambda}
       (\pi_\alpha^\lambda)^{-1}(E_{n,\alpha})
       \quad(0<\lambda\le\omega_1\text{ a limit}).
       \label{eq:E-limit-definition}
\end{align}
In \eqref{eq:E-successor}, $z_n^+(\alpha)$ is the new $\C$-block
of the $n$th auxiliary output. The zero-slice assumption
\eqref{eq:zero-slices} and induction give
\begin{equation}\label{eq:K-in-E}
 K_{n,\alpha}\subseteq E_{n,\alpha}.
\end{equation}
Also,
\begin{equation}\label{eq:E-coherence}
 \pi_\alpha^\delta(E_{n,\delta})\subseteq E_{n,\alpha}
 \qquad(\alpha\le\delta\le\omega_1).
\end{equation}
For successive indices, the left side is contained in
$K_{n,\alpha}$ by definition, and then in $E_{n,\alpha}$ by
\eqref{eq:K-in-E}. The limit step follows from the defining
intersection; transfinite induction gives the general case.

For every nonzero limit $\lambda\le\omega_1$ we have
\begin{equation}\label{eq:E-equals-K-limit}
 E_{n,\lambda}=K_{n,\lambda}.
\end{equation}
The inclusion from right to left follows from
\eqref{eq:K-coherence} and \eqref{eq:K-in-E}. For the other
inclusion, a point of $E_{n,\lambda}$ has its
$(\alpha+1)$-prefix in $E_{n,\alpha+1}$ for every
$\alpha<\lambda$, since $\lambda$ is a limit. Its
$\alpha$-prefix therefore lies in $K_{n,\alpha}$. Now use
\eqref{eq:K-limit}. We do not assert equality of $E_{n,\alpha}$
and $K_{n,\alpha}$ at arbitrary successor stages. Notice that
\eqref{eq:E-coherence} implies that once an output prefix leaves
its exceptional set, no longer prefix can re-enter.

\medskip\noindent\emph{The exact fibre invariant.}
Put
\[
 X_\alpha=\B\times P_\alpha\times(P_\alpha)^{\N},
 \qquad P_\alpha=\B^\alpha,
\]
and let $\rho_\alpha^\delta:X_\delta\to X_\alpha$ preserve
the first factor and restrict all remaining tuples. We construct
continuous maps $H_\alpha:X_\alpha\to Y_\alpha$ such that
\begin{equation}\label{eq:H-compatibility}
 \pi_\alpha^\delta H_\delta=H_\alpha\rho_\alpha^\delta
 \qquad(\alpha\le\delta\le\omega_1).
\end{equation}
Let $\operatorname{pr}_{V\times\C^\alpha}$ denote the projection
of $Y_\alpha$ onto its first two factors. We require
\begin{equation}\label{eq:primary-outputs}
 \operatorname{pr}_{V\times\C^\alpha}
 \bigl(H_\alpha(t,s,(x_n)_n)\bigr)
 =(b(t),c^\alpha(s)).
\end{equation}
For every \emph{countable} $\alpha$ and every
$y=(v,a,(z_n)_n)\in Y_\alpha$, the invariant is the actual set
equality
\begin{equation}\label{eq:fibre-product}
 H_\alpha^{-1}(y)=
 \{b^{-1}(v)\}\times\{(c^\alpha)^{-1}(a)\}
       \times\prod_{n\in\N}L_{n,\alpha}(y),
\end{equation}
where
\begin{equation}\label{eq:L-definition}
 L_{n,\alpha}(y)=
 \begin{cases}
 Z_\alpha,&y\in E_{n,\alpha},\\
 \{\ell_{n,\alpha}(y)\},&y\notin E_{n,\alpha}
 \end{cases}
\end{equation}
for some $\ell_{n,\alpha}(y)\in P_\alpha$ in the second case.
The inverse symbols in \eqref{eq:fibre-product} specify unique
points as sets; no continuity of those inverses or of the
functions $\ell_{n,\alpha}$ is being asserted. Since $Z_\alpha$
is nonempty at each countable stage, this invariant implies that
$H_\alpha$ is onto at every such stage. This nonemptiness will
be used in the countable limit step.

\medskip\noindent\emph{Initial stage.}
After suppressing singleton factors, let $H_0=b$. Here
$E_{n,0}=V$ and $Z_0=P_0$ is a singleton, so all the required
conditions hold. The choice $E_{n,0}=V$ is intentional, even if
$K_{n,0}$ is a proper subset of $V$.

\medskip\noindent\emph{The two extensions at a successor stage.}
Fix $\alpha<\omega_1$ and $n\in\N$. Abbreviate
$E=E_{n,\alpha}$ and $K=K_{n,\alpha}$. These are closed
zero-dimensional Polish subspaces of $Y_\alpha$, and $K\subseteq E$.
Recall that $W_\alpha=Z_\alpha\times\B\cong\B$.

The two extensions address different types of old fibres. On $E$ an old auxiliary tuple ranges over the entire $Z_\alpha$, so the old tuple and the new Baire variable must be treated together: on $K$ the zero value retains precisely the next exceptional set, while on $E\setminus K$ every new value must determine a unique pair. Outside $E$ the old tuple is already unique, so only the new Baire variable has to be determined. A single extension with bijectivity in that new variable would not resolve the entire old exceptional factor on $E\setminus K$.

On $K\times W_\alpha$ prescribe
$(y,x,u)\mapsto r_\alpha(x,u)$, where
$r_\alpha$ is given by Lemma~\ref{lem:exceptional-collapse}.
Apply Lemma~\ref{lem:parameter-extension} with base $E$, closed
subset $K$, and second factor $W_\alpha$, using its homeomorphism
with $\B$. We obtain a continuous map
\[
 R_{n,\alpha}:E\times Z_\alpha\times\B\longrightarrow\C
\]
such that it equals $r_\alpha(x,u)$ for $y\in K$ and, for each
$y\in E\setminus K$, the map
$(x,u)\mapsto R_{n,\alpha}(y,x,u)$ is a condensation
$Z_\alpha\times\B\to\C$. It is the \emph{pair} $(x,u)$
that is used in this bijectivity assertion.

Now apply Lemma~\ref{lem:parameter-extension} with base
$Y_\alpha\times P_\alpha$, closed subset $E\times Z_\alpha$,
and second factor $\B$. It extends $R_{n,\alpha}$ to a continuous
map
\[
 q_{n,\alpha}:Y_\alpha\times P_\alpha\times\B\longrightarrow\C
\]
such that $u\mapsto q_{n,\alpha}(y,x,u)$ is a condensation whenever
$(y,x)\notin E\times Z_\alpha$. Here bijectivity is in the new
variable $u$ alone. Both bases are zero-dimensional Polish spaces because $\alpha$ is countable; the prescribed subsets are closed in the corresponding bases. Empty $K$ or $E$ is covered by the empty-domain cases
of the same lemma.

\medskip\noindent\emph{The successor map and its whole fibre.}
Write an input of $X_{\alpha+1}$ as
\[
 (t,(s,w),((x_n,u_n))_{n\in\N}),
 \quad s,x_n\in P_\alpha,\quad w,u_n\in\B,
\]
and write $H_\alpha(t,s,(x_n)_n)=y=(v,a,(z_n)_n)$. Define
\begin{equation}\label{eq:H-successor}
 \begin{split}
 &H_{\alpha+1}(t,(s,w),((x_n,u_n))_n)\\
 &\qquad=
 \bigl(v,(a,c(w)),
          ((z_n,q_{n,\alpha}(y,x_n,u_n)))_{n\in\N}\bigr).
 \end{split}
\end{equation}
Continuity follows coordinatewise from continuity of $H_\alpha$,
$c$ and the $q_{n,\alpha}$. Compatibility and the prescribed first
two outputs follow from the formula.

Fix a desired output
$y^+=(v,a^+,(z_n^+)_n)\in Y_{\alpha+1}$ and let
$y=\pi_\alpha^{\alpha+1}(y^+)$. The old first two inputs $t,s$
are uniquely fixed by \eqref{eq:fibre-product}. The new primary
input $w$ is uniquely fixed by $c(w)=a^+(\alpha)$. Conditional on
the old output $y$, the old auxiliary inputs form exactly the
product in \eqref{eq:fibre-product}. Every new equation for index
$n$ uses this fixed old \emph{output} $y$, together with only
$x_n$ and $u_n$. It therefore imposes no relation between the
auxiliary inputs for different indices. More explicitly, write
\[
 S_n(y^+)=\{(x,u)\in L_{n,\alpha}(y)\times\B:
 q_{n,\alpha}(y,x,u)=z_n^+(\alpha)\}.
\]
Under the natural grouping of coordinates, the full new fibre has the uniquely determined first factors $t,s,w$ and the product of the sets $S_n(y^+)$. Indeed, a preimage must satisfy each displayed equation and its old prefix belongs to the old fibre by compatibility. Conversely, any family of pairs from these sets, with $t,s,w$ fixed, realizes the old output $y$ by the full old fibre identity. Thus the parameter in each new equation is exactly $y$, and every new output coordinate is the prescribed one.
For each $n$ we have the following four cases.

\begin{enumerate}
\item If $y\notin E_{n,\alpha}$, then $x_n$ is already unique.
The point $(y,x_n)$ lies outside $E_{n,\alpha}\times Z_\alpha$,
so bijectivity of $q_{n,\alpha}(y,x_n,\cdot)$ gives a unique
$u_n$ for the new output. Since $K_{n,\alpha}\subseteq
E_{n,\alpha}$, the new output is outside $E_{n,\alpha+1}$.

\item If $y\in E_{n,\alpha}\setminus K_{n,\alpha}$, then
$x_n$ ranges over $Z_\alpha$. The first extension is bijective
in $(x_n,u_n)$ on $Z_\alpha\times\B$, so every prescribed
new output gives exactly one such pair. Again the new output
is outside $E_{n,\alpha+1}$.

\item If $y\in K_{n,\alpha}$ and $z_n^+(\alpha)=0$, then
$x_n$ ranges over $Z_\alpha$ and the new output map is $r_\alpha$.
Its zero fibre is exactly $Z_{\alpha+1}$, and $y^+\in E_{n,\alpha+1}$.

\item If $y\in K_{n,\alpha}$ and $z_n^+(\alpha)\ne0$, the same map
$r_\alpha$ gives exactly one pair $(x_n,u_n)$. In this case
$y^+\notin E_{n,\alpha+1}$, by \eqref{eq:E-successor}.
\end{enumerate}

These independent equations prove both inclusions in
\eqref{eq:fibre-product} at $\alpha+1$: every preimage has the
specified factors, and every combination of allowed factors
realizes all the output equations. In particular, this is a
statement about the whole fibre and not merely about its
coordinate projections. Its factors are nonempty, so
$H_{\alpha+1}$ is onto.

\medskip\noindent\emph{A countable limit stage.}
Let $0<\lambda<\omega_1$ be a limit. Define $H_\lambda$ by
compatibility. Every output block of index $\xi<\lambda$ is
provided by $H_{\xi+1}$ composed with the continuous prefix
projection. The $V$-coordinate is $b(t)$. Thus $H_\lambda$ is
continuous, satisfies \eqref{eq:primary-outputs}, and is compatible
with all preceding maps.

For $y\in Y_\lambda$, abbreviate
$y_\alpha=\pi_\alpha^\lambda(y)$. We have, directly from the
definition,
\begin{equation}\label{eq:fibre-limit}
 H_\lambda^{-1}(y)=
 \bigcap_{\alpha<\lambda}
       (\rho_\alpha^\lambda)^{-1}(H_\alpha^{-1}(y_\alpha)).
\end{equation}
Substituting the known identities \eqref{eq:fibre-product}, the
first two factors are fixed, and the $n$th auxiliary factor is
exactly
\begin{equation}\label{eq:M-factor}
 M_n=\{x\in P_\lambda:x\res\alpha\in
                L_{n,\alpha}(y_\alpha)\text{ for all }\alpha<\lambda\}.
\end{equation}
This factorization just interchanges the universal quantifiers
over $\alpha$ and $n$ in a set of independent coordinate
conditions. It assumes no theorem about surjectivity of inverse
systems.

If $y\in E_{n,\lambda}$, every $y_\alpha$ belongs to
$E_{n,\alpha}$, so \eqref{eq:M-factor} equals $Z_\lambda$ by
\eqref{eq:Z-limit}. Otherwise choose $\beta<\lambda$ such that
$y_\beta\notin E_{n,\beta}$. By \eqref{eq:E-coherence},
$y_\gamma\notin E_{n,\gamma}$ for
$\beta\le\gamma<\lambda$. Write the corresponding singleton
as $L_{n,\gamma}(y_\gamma)=\{\ell_\gamma\}$.

These singleton tuples are compatible. Indeed, if
$\beta\le\gamma<\delta<\lambda$, the full fibre
$H_\delta^{-1}(y_\delta)$ is nonempty by the previously proved
countable-stage invariant. Choose a point of it. Its prefix lies
in $H_\gamma^{-1}(y_\gamma)$ by compatibility, and therefore
$\ell_\delta\res\gamma=\ell_\gamma$. The same argument shows
that their earlier restrictions belong to the factors
$L_{n,\alpha}(y_\alpha)$ for $\alpha<\beta$. Hence
$\bigcup_{\beta\le\gamma<\lambda}\ell_\gamma$ belongs to
$M_n$, and is its only member because those lengths are cofinal
in $\lambda$.

We have proved \eqref{eq:fibre-product} at $\lambda$.
All its factors are nonempty, including $Z_\lambda$, so the
map is onto. This completes the induction through all countable
stages. In particular, we never extended an arbitrary element of
$Z_\alpha$ to $Z_{\alpha+1}$.

\medskip\noindent\emph{The final stage.}
\emph{Continuity.} Define $H_{\omega_1}$ by compatibility as above. Each output coordinate is obtained from a continuous map at a countable stage, so $H_{\omega_1}$ is continuous.
For $y\in Y$, write $y_\alpha=\pi_\alpha(y)$.
\emph{Empty fibres.} If $y\in K_n$ for some $n$, every countable prefix lies in
$K_{n,\alpha}\subseteq E_{n,\alpha}$. A preimage of $y$ would
have its $n$th auxiliary tuple in $Z_\alpha$ at every countable
prefix, and hence in $Z_{\omega_1}$ by \eqref{eq:Z-limit}.
This is impossible by \eqref{eq:Z-empty}.

\emph{Singleton fibres.} Suppose instead that $y\notin\bigcup_nK_n$. Since
$K_n=E_{n,\omega_1}$, for every $n$ there is a countable
$\beta_n$ with $y_{\beta_n}\notin E_{n,\beta_n}$. Put
\begin{equation}\label{eq:common-stage}
 \beta=\sup_{n\in\N}\beta_n<\omega_1.
\end{equation}
For all $\beta\le\gamma<\omega_1$, coherence implies that
$y_\gamma$ is outside every $E_{n,\gamma}$. The complete
countable-stage fibre is consequently a singleton; denote its
point by $p_\gamma$. Compatibility gives
$\rho_\gamma^\delta(p_\delta)=p_\gamma$ for
$\beta\le\gamma<\delta<\omega_1$. Write $p_\gamma=(t_\gamma,s_\gamma,(x_{n,\gamma})_n)$. Compatibility gives a common value $t$ of all $t_\gamma$. Define
\[
 p=\left(t,\bigcup_{\beta\le\gamma<\omega_1}s_\gamma,
        \left(\bigcup_{\beta\le\gamma<\omega_1}x_{n,\gamma}\right)_{n\in\N}\right).
\]
The cofinal compatible tuples define all coordinates of $p\in X_{\omega_1}$.
Every output prefix of this element agrees with the corresponding
prefix of $y$, so it maps to $y$. Any other preimage would have
the same prefixes and must be equal to it. We have proved
\begin{equation}\label{eq:final-fibres}
 |H_{\omega_1}^{-1}(y)|=
 \begin{cases}
 0,&y\in\bigcup_{n\in\N}K_n,\\
 1,&y\notin\bigcup_{n\in\N}K_n.
 \end{cases}
\end{equation}

Finally,
\[
 X_{\omega_1}=
 \B\times\B^{\omega_1}\times(\B^{\omega_1})^{\N}
 \cong\N^{\omega_1}
\]
by a coordinate bijection: the total number of discrete
$\N$-coordinates is $\aleph_1$. A continuous map into $Y$,
regarded as a map into its image with the subspace topology, is
still continuous. Thus \eqref{eq:final-fibres} gives the required
condensation onto $Y\setminus\bigcup_nK_n$.
\end{proof}

\section{Completion and immediate consequences}\label{sec:completion}

\begin{proof}[Proof of Theorem~\ref{thm:main-gdelta}]
Write $B=\bigcap_{n\in\N}U_n\ne\varnothing$ in
$2^{\omega_1}$. Proposition~\ref{prop:polish-parameter-reduction}
gives a countable set $A$ and
$V=\pi_A(B)=V_0\times2^L$, with $V_0$ nonempty
zero-dimensional Polish and $L$ countably infinite. Every section
$(U_n)_a$, $a\in V$, is $G_\delta$-dense. By
Lemma~\ref{lem:parameter-absorption}, there is a condensation
$b:\B\to V$.

Put $X=V\times2^{\omega_1\setminus A}$ and
$F_n=X\setminus U_n$. Since $\pi_A(B)=V$, the set $B$ is
contained in $X$ and
$B=X\setminus\bigcup_nF_n$. Proposition~\ref{prop:simultaneous-normalization}
provides a homeomorphism
\[
 h:X\longrightarrow
 Y=V\times\C^{\omega_1}\times(\C^{\omega_1})^{\N}
\]
such that $K_n=h(F_n)$ is closed and satisfies
\eqref{eq:zero-slices}. Thus
$h(B)=Y\setminus\bigcup_nK_n$. Apply
Proposition~\ref{prop:countable-exceptions}, and compose the
resulting condensation with $h^{-1}|_{h(B)}$.
\end{proof}

\begin{corollary}\label{cor:complements}
If $(F_n)_{n\in\N}$ is a sequence of closed subsets of
$2^{\omega_1}$ whose union is proper, then
$\N^{\omega_1}$ condenses onto
$2^{\omega_1}\setminus\bigcup_nF_n$. In particular, it
condenses onto every nonempty open subset and onto the complement
of every countable subset of $2^{\omega_1}$.
\end{corollary}
\begin{proof}
The indicated complement is a nonempty countable intersection of
open sets. For the last assertion, singletons are closed and a
countable set cannot exhaust $2^{\omega_1}$.
\end{proof}

\begin{corollary}\label{cor:other-sources}
Let $M$ be a noncompact zero-dimensional Polish space. Then $M^\kappa\cong\N^\kappa$ for every infinite cardinal $\kappa$. In particular, at the first uncountable cardinal, $M^{\omega_1}$ condenses onto every nonempty $G_\delta$-subset of $2^{\omega_1}$.
\end{corollary}
\begin{proof}
The space $M^\N$ is nonempty zero-dimensional Polish and nowhere
locally compact. Every nonempty basic open set contains a closed
copy of $M$ obtained by fixing all but one unrestricted
coordinate. Such a copy is noncompact. Hence
$M^\N\cong\B$ by Lemma~\ref{lem:characterizations}(2).
Reindexing gives
$M^\kappa\cong(M^\N)^\kappa\cong\B^\kappa\cong\N^\kappa$.
Apply Theorem~\ref{thm:main-gdelta}.
\end{proof}

\section{A corollary for non-Archimedean function spaces}
The following consequence is obtained by topological identifications of the source and target with the products in Theorem~\ref{thm:main-gdelta}; it is not a separate condensation construction. For a discrete space $D$ and a valued field $L$, write $C_p(D,L)$ for the space of all continuous $L$-valued functions on $D$ with the topology of pointwise convergence. Thus $C_p(D,L)=L^D$. The conclusion gives condensations from this pointwise function space onto prescribed nonempty $G_\delta$-subspaces of the weak$^*$ dual ball, not merely onto the whole ball; no linearity of the condensation is asserted.

\begin{corollary}\label{cor:local}
Let $L$ be a nondiscrete, locally compact field with a nontrivial non-Archimedean absolute value, endowed with the topology induced by that absolute value, and let $D$ be discrete of cardinality $\aleph_1$. Put $E=c_0(D,L)$, with the supremum norm, and let $B_{E'}$ be the closed unit ball of its continuous dual. Then $C_p(D,L)$ condenses onto every nonempty $G_\delta$-subset of $(B_{E'},w^*)$.
\end{corollary}
\begin{proof}
Write $B_L=\{a\in L:|a|\le1\}$. Choose a compact neighborhood of zero and a closed ball $aB_L$ contained in it, where $a\ne0$. Such a ball exists by nontriviality of the absolute value. The ball $aB_L$ is closed and therefore compact, and multiplication by $a^{-1}$ shows that $B_L$ is compact. The field is metrizable and $B_L$ is zero-dimensional with no isolated points, so $B_L\cong\C$ by Lemma~\ref{lem:characterizations}(1). For $0<|a|<1$ we have $L=\bigcup_{n\ge0}a^{-n}B_L$. Each of these compact metrizable balls is separable, so $L$ is separable. The field is noncompact because $|a^{-n}|\to\infty$. It is complete: a Cauchy sequence is eventually in a compact translate of $B_L$ and thus has a convergent subsequence, which forces the sequence itself to converge. As ultrametric balls give a clopen base, $L$ is a noncompact zero-dimensional Polish space.

Each $\varphi\in B_{E'}$ determines a tuple $(\varphi(e_d))_{d\in D}\in B_L^D$, where $e_d$ is the unit vector. Conversely, for $u=(u_d)\in B_L^D$, the formula
\[
 \varphi_u(x)=\sum_{d\in D}u_dx_d\qquad(x\in c_0(D,L))
\]
defines a continuous functional of norm at most one. The sum is the limit of the net of finite partial sums; this limit exists because for each $\varepsilon>0$ only finitely many $|x_d|$ are at least $\varepsilon$, and the field is complete. These correspondences are inverse, by density of the finitely supported vectors. In fact $\|\varphi_u\|=\sup_{d\in D}|u_d|$: the upper bound follows from the ultrametric inequality and the lower bound by evaluating at the unit vectors $e_d$. The coordinate evaluations are weak$^*$-continuous. Conversely, for every $x\in E$ the finite partial sums converge uniformly on $B_L^D$, by the ultrametric tail estimate $|\sum_{d\notin F}u_dx_d|\le\sup_{d\notin F}|x_d|$ for finite $F\subseteq D$. Hence $u\mapsto\varphi_u(x)$ is product-continuous. This identifies $(B_{E'},w^*)$ homeomorphically with $B_L^D\cong2^{\omega_1}$. The source is homeomorphic to $\N^{\omega_1}$ by Corollary~\ref{cor:other-sources}; apply Theorem~\ref{thm:main-gdelta}.
\end{proof}

\section{Larger cardinals}
The main theorem settles all nonempty $G_\delta$ targets at weight $\aleph_1$. The results below distinguish three issues: condensations for targets determined by few coordinates (Proposition~\ref{prop:gdelta} and Corollary~\ref{cor:smallcoordinates}), surjections and their fibres for punctured cubes (Proposition~\ref{prop:extent} and Theorem~\ref{thm:singleton}), and obstructions to the proof method or to the existence of the maps themselves. None gives a general condensation theorem for $\kappa>\aleph_1$. For comparison, complements of proper closed $G_\delta$-subsets admit an elementary positive result at every uncountable weight; these complements need not be $G_\delta$-dense.

\begin{proposition}\label{prop:gdelta}
Let $\kappa$ be uncountable and $F$ a proper closed $G_\delta$-subset of $2^\kappa$. Then there is a condensation $\N^\kappa\to2^\kappa\setminus F$.
\end{proposition}
\begin{proof}
If $F=\varnothing$ take a power of Lemma~\ref{lem:basic-condensation}. Otherwise write $F=\bigcap_n U_n$ with $U_n$ open. Compactness and the clopen base yield clopen sets $V_n$ with $F\subseteq V_n\subseteq U_n$, each determined by finitely many coordinates. Thus $F$ is determined by a countable set of coordinates $A$, which we may enlarge to be countably infinite. In the decomposition $2^\kappa=2^A\times2^{\kappa\setminus A}$ the complement is $V\times 2^{\kappa\setminus A}$ with $V$ a nonempty open subset of $2^A$. Disjointifying a countable cover by basic clopen sets partitions $V$ into a nonempty finite or countably infinite family of clopen Cantor spaces. The complement is accordingly a topological sum of that many copies of $2^\kappa$. On the other hand, $\N^\kappa$ is a topological sum of any prescribed nonzero finite or countably infinite number of copies of itself, by partitioning one discrete factor into that many infinite sets. Combine the corresponding condensations obtained from Lemma~\ref{lem:basic-condensation} by taking powers.
\end{proof}

\begin{corollary}\label{cor:smallcoordinates}
Let $\kappa\ge\aleph_1$. Suppose that a nonempty $G_\delta$-subset $B$ of $2^\kappa$ is determined by at most $\aleph_1$ coordinates. Then there is a condensation $\N^\kappa\to B$. In particular, this holds when $B$ is a nonempty intersection of countably many open sets, each of which is a union of at most $\aleph_1$ basic cylinders.
\end{corollary}
\begin{proof}
If $\kappa=\aleph_1$, apply Theorem~\ref{thm:main-gdelta}. Otherwise choose $A\subseteq\kappa$ of cardinality $\aleph_1$ containing the determining coordinates. Then $B=V\times2^{\kappa\setminus A}$ for a nonempty subset $V$ of $2^A$. Fixing any point in the second factor embeds $2^A$ continuously into the full cube; its inverse image of $B$ is $V$, so $V$ is $G_\delta$. Apply Theorem~\ref{thm:main-gdelta} on $A$, and take the product with the elementary condensation $\N^{\kappa\setminus A}\to2^{\kappa\setminus A}$, obtained by reindexing powers of Lemma~\ref{lem:basic-condensation}. The last assertion follows by collecting the finite coordinate sets of the cylinders used in the countable intersection; their union has cardinality at most $\aleph_1$.
\end{proof}

The following observation gives a necessary condition at higher cardinals and isolates the difference between surjectivity and bijectivity.

\begin{proposition}\label{prop:extent}
For every infinite cardinal $\kappa$, the following are equivalent:
\begin{enumerate}
 \item There is a continuous surjection $\N^\kappa\to2^\kappa\setminus\{0\}$.
 \item The product $\N^\kappa$ contains a closed discrete subspace of cardinality $\kappa$.
\end{enumerate}
\end{proposition}
\begin{proof}
For the forward implication, the family of unit vectors $(e_i)_{i<\kappa}$ is closed and discrete in the punctured cube. Choose a preimage of each vector. These chosen points are discrete and closed: inside the closed inverse image of the family, the individual vector fibres are clopen, and the chosen singleton is closed in each fibre.

Conversely, let $D=\{d_i:i<\kappa\}$ be closed discrete in $P=\N^\kappa$. Choose a clopen neighborhood $V_i$ of $d_i$ with $V_i\cap D=\{d_i\}$. The clopen cylinder base of $P$ has cardinality $\kappa$. Thus $P\setminus D$ is covered by at most $\kappa$ clopen sets disjoint from $D$. Assign each of these sets to a distinct index $i<\kappa$ and let $U_i$ be $V_i$ united with the assigned set, if any. Then $(U_i)$ covers $P$, each $U_i$ is clopen, and $U_i\cap D=\{d_i\}$. The map $s:P\to2^\kappa$ given by $s(x)(i)=\mathbf1_{U_i}(x)$ is continuous, never zero and satisfies $s(d_i)=e_i$. Consequently
\[
 (x,z)\longmapsto (\max\{s(x)(i),z(i)\})_{i<\kappa}
\]
is a continuous surjection $P\times2^\kappa\to2^\kappa\setminus\{0\}$. Indeed, for a prescribed nonzero $z$ choose an index $i$ with $z(i)=1$ and take $x=d_i$. There is a continuous surjection $\N^\kappa\to2^\kappa$ by taking a coordinatewise surjection $\N\to 2$, and $P\times P\cong P$. This proves the result.
\end{proof}

The closed-discrete-set condition in Proposition~\ref{prop:extent} holds in particular at $\aleph_2$. Mycielski's theorem \cite{Mycielski}, as stated in the review MR0211871, asserts that $\N^\kappa$ contains a subset of cardinality $\kappa$ without accumulation points for every infinite cardinal $\kappa$ below the first uncountable weakly inaccessible cardinal. Such a set is closed and discrete. The original article has not been consulted directly. See also Gorelic \cite[p.~613]{Gorelic}, Basile \cite[proof of Corollary 4.14]{Basile} and Lipparini \cite[pp.~239, 245]{Lipparini}. We include a short direct proof of this classical closed-discrete-set assertion, including the singular-cardinal step, for self-containedness. No new range of cardinals is claimed here; the consequence for surjectivity follows by applying Proposition~\ref{prop:extent}, and Theorem~\ref{thm:singleton} below supplies the additional control of fibres.

\begin{lemma}\label{lem:successorextent}
If $\mu$ is an infinite cardinal and $D_\mu$ is discrete of cardinality $\mu$, then $D_\mu^{\mu^+}$ contains a closed discrete subset of cardinality $\mu^+$.
\end{lemma}
\begin{proof}
Use the homeomorphic product $(2\times\mu)^{\mu^+}$, where $\mu$ is discrete. For each $\beta<\mu^+$ choose an injection $j_\beta:\beta\to\mu$. Define $x_\beta$ by
\[
 x_\beta(\alpha)=
 \begin{cases}
 (0,j_\beta(\alpha)),&\alpha<\beta,\\
 (1,j_{\alpha+1}(\beta)),&\beta\le\alpha<\mu^+.
 \end{cases}
\]
The expression $j_{\alpha+1}(\beta)$ is defined since $\beta<\alpha+1<\mu^+$. Put $D=\{x_\beta:\beta<\mu^+\}$. If $x\in\overline D$ has some coordinate with first entry $1$, say $x(\alpha)=(1,\xi)$, then the clopen cylinder specifying this coordinate meets $D$ in exactly one point, by injectivity of $j_{\alpha+1}$ and the closure assumption. Since a singleton is closed, $x$ must be that point. Otherwise $x(\alpha)=(0,\xi_\alpha)$ for all $\alpha<\mu^+$. For distinct $\alpha,\gamma$, a two-coordinate cylinder around $x$ meets $D$ in an $x_\beta$ with $\beta>\max(\alpha,\gamma)$. Thus $\xi_\alpha\ne\xi_\gamma$ by injectivity of $j_\beta$, giving an impossible injection $\mu^+\to\mu$. Consequently $D$ is closed. Each $x_\beta$ is isolated in $D$ by its $\beta$th coordinate, which also shows that these points are all distinct.
\end{proof}

\begin{corollary}\label{cor:mycielski}
Suppose that $\kappa$ is infinite and that there is no uncountable regular limit cardinal at most $\kappa$. Then $\N^\kappa$ contains a closed discrete subset of cardinality $\kappa$, and hence admits a continuous surjection onto $2^\kappa\setminus\{0\}$. In particular this holds for every $\kappa=\aleph_n$ with $n<\omega$.
\end{corollary}
\begin{proof}
For $\kappa=\aleph_0$, use a single unrestricted discrete coordinate and fix the rest. Induct on the infinite cardinals in the indicated range. For $\kappa=\mu^+$, a closed discrete subspace of size $\mu$ in $\N^\mu$ is a closed copy of $D_\mu$. Lemma~\ref{lem:successorextent} gives a closed discrete subspace of size $\mu^+$ in $D_\mu^{\mu^+}$, which embeds as a closed subspace of $(\N^\mu)^{\mu^+}\cong\N^{\mu^+}$.

If $\kappa$ is an uncountable limit cardinal in the indicated range, it is singular. Set $\delta=\operatorname{cf}(\kappa)<\kappa$ and choose a cofinal family of infinite cardinals $\kappa_i<\kappa$ for $i<\delta$. By induction choose a closed discrete set $\{d_i:i<\delta\}$ in $\N^\delta$ and closed discrete subsets $E_i$ of $\N^\kappa$ of size $\kappa_i$ (pad the corresponding products with fixed coordinates). The subspace $\bigcup_{i<\delta}\{d_i\}\times E_i$ is closed and discrete in $\N^\delta\times\N^\kappa$: it is closed inside the closed product $\{d_i:i<\delta\}\times\N^\kappa$, whose first factor is discrete. The product $\{d_i:i<\delta\}\times\N^\kappa$ is a topological sum of its clopen slices. The constructed union has cardinality $\sum_{i<\delta}\kappa_i=\kappa$ by cofinality and $\delta\cdot\kappa=\kappa$. Thus the required cardinality is attained, not merely approached by smaller closed discrete subsets. Reindex the ambient product and apply Proposition~\ref{prop:extent}.
\end{proof}

We can strengthen the surjection statement by specifying a large set of singleton fibres. For $y\in2^\kappa$, let $\operatorname{supp}(y)=\{i<\kappa:y(i)=1\}$.

\begin{theorem}\label{thm:singleton}
Let $\kappa$ be an infinite cardinal such that $\N^\kappa$ contains a closed discrete subset of cardinality $\kappa$. There is a continuous surjection $F:\N^\kappa\to2^\kappa\setminus\{0\}$ satisfying
\[
 |F^{-1}(y)|=1\qquad\text{whenever }0<|\operatorname{supp}(y)|<\kappa.
\]
\end{theorem}
\begin{proof}
Write $P=\N^\kappa$ and fix a closed discrete set $\{d_i:i<\kappa\}$. We may assume, after replacing this set by a graph in $P\times2^\kappa\cong P$, that for every $\gamma<\kappa$ there is a clopen set $T_\gamma\subset P$ with trace $\{d_i:i<\gamma\}$ on $D=\{d_i:i<\kappa\}$. Indeed, take the graph of $d_i\mapsto (\mathbf{1}_{i<\gamma})_{\gamma<\kappa}$ over the closed discrete set. The graph is closed in $D\times2^\kappa$ because the map from the discrete set $D$ is continuous and $2^\kappa$ is Hausdorff. As $D$ is closed in $P$, the graph is closed in $P\times2^\kappa$; it is discrete by the isolating neighborhoods in $D$. The clopen sets in the added factor give the required traces. Set also $T_\kappa=P$.

Choose clopen $V_i\subset P$ with $V_i\cap D=\{d_i\}$. Cover $P\setminus D$ by at most $\kappa$ clopen cylinders disjoint from $D$. Assign each cylinder in the cover a set of $\kappa$ indices below $\kappa$, pairwise disjoint for distinct cylinders; this is possible by $\kappa\cdot\kappa=\kappa$. For each $i$, form $U_i$ by adjoining to $V_i$ the cylinder assigned to $i$, if there is one. Define $s:P\to2^\kappa$ by $s(x)(i)=\mathbf1_{U_i}(x)$. It is continuous, with $s(d_i)=e_i$ and $|\operatorname{supp}s(x)|=\kappa$ for $x\notin D$.

For $i<j<\kappa$, let $k(i,j)=\operatorname{otp}([i+1,j))$, where $\operatorname{otp}$ denotes order type and $[i+1,j)=\{\xi:i+1\le\xi<j\}$ is a half-open ordinal interval. Then $(i+1)+k(i,j)=j$, and strict monotonicity of ordinal addition in its right argument gives uniqueness. For fixed $i$, the map $k\mapsto(i+1)+k$ is an order isomorphism of $\kappa$ onto $[i+1,\kappa)$. Indeed, the latter interval has cardinality $\kappa$ and order type at most $\kappa$, since it is a subset of $\kappa$. As $\kappa$ is an initial ordinal, this order type must equal $\kappa$. For each fixed $j$ the set
\[
 L_j=\{k(i,j):i<j\}
\]
is finite. For $j=0$ this is immediate. For $j>0$, write its Cantor normal form as $j=\omega^{\eta_1}c_1+\cdots+\omega^{\eta_r}c_r$, where $\eta_1>\cdots>\eta_r$ and $0<c_\ell<\omega$. Suppose $j=a+k$ with $k>0$, and write $k=\omega^\zeta c+v$ with $0<c<\omega$ and $v<\omega^\zeta$. The addition rule deletes the terms of $a$ of exponent less than $\zeta$, adds the coefficient at $\zeta$ (if present) to $c$, and preserves all terms of greater exponent. Thus $\zeta=\eta_\ell$ for some $\ell$, the resulting coefficient at $\zeta$ is $c_\ell\ge c$, and the terms after $\zeta$ are precisely those of $v$. Consequently
\[
 k=\omega^{\eta_\ell}c+\omega^{\eta_{\ell+1}}c_{\ell+1}+\cdots+\omega^{\eta_r}c_r,
 \qquad 1\le c\le c_\ell,
\]
where the sum after the first term is empty if $\ell=r$. There are at most $\sum_{\ell=1}^r c_\ell$ positive possibilities and also the possibility $k=0$. This proves finiteness of $L_j$. For each $k\in L_j$ the set $\{i<j:k(i,j)=k\}$ is an ordinal interval, because $i\mapsto(i+1)+k$ is nondecreasing. It therefore has the form $[a,b)=\{i:a\le i<b\}$ with $a$ its least element and $b=\sup\{i+1:i<j,\ k(i,j)=k\}\le j$. Choose the clopen set $W_{j,k}=T_b\setminus T_a$; its trace on $D$ consists exactly of the corresponding $d_i$.

For $(x,z)\in P\times2^\kappa$, define
\[
 H(x,z)(j)=s(x)(j)\vee
 \bigvee_{k\in L_j}\bigl(\mathbf1_{W_{j,k}}(x)\wedge z(k)\bigr),
\]
where the empty disjunction is zero and $\vee,\wedge$ denote the Boolean operations. Each coordinate is continuous by finiteness of $L_j$. Also $H(x,z)\ge s(x)$ coordinatewise, so the value is never zero. On the fibre over $d_i$ the formula reduces to
\[
 H(d_i,z)(j)=
 \begin{cases}0,&j<i,\\1,&j=i,\\z(k(i,j)),&j>i.\end{cases}
\]
Thus $H(d_i,\cdot)$ is a bijection onto the set of points whose first nonzero coordinate is $i$. These sets partition the punctured cube, so $H$ is onto and its restriction to $D\times2^\kappa$ is bijective. But every $H(x,z)$ with $x\notin D$ has support of cardinality $\kappa$. This proves the claimed singleton-fibre property. Finally compose with the inverse of the homeomorphism $P\times2^\kappa=\N^\kappa\times2^\kappa\cong(\N\times2)^\kappa\cong\N^\kappa$, where the last map is induced by the bijection $(n,\varepsilon)\mapsto2n-1+\varepsilon$ between discrete spaces. This gives $F$.
\end{proof}

\begin{remark}
Theorem~\ref{thm:singleton} does not establish injectivity over points of support of cardinality $\kappa$. In particular, it does not settle the condensation problem at $\aleph_2$. The restriction of $H$ to $D\times2^\kappa$ is already surjective, so changing only the values outside that subspace cannot make $H$ a bijection.
\end{remark}

\begin{proposition}\label{prop:ultrafilter}
Suppose there exists a free countably complete ultrafilter on a set $J$. For every set $I$ there is no continuous surjection $\N^I\to2^J\setminus\{0\}$.
\end{proposition}
\begin{proof}
Let $\mathcal U$ be such an ultrafilter and suppose $h$ is a continuous surjection. Choose $x_j\in\N^I$ with $h(x_j)=e_j$. Each coordinate $i\in I$ is constant on a member of $\mathcal U$: the countable partition of $J$ into the fibres of $j\mapsto x_j(i)$ has a member in $\mathcal U$ by countable completeness. Let $x(i)$ be the corresponding constant value. Finite intersections show that $x_j\to x$ along $\mathcal U$ in the product topology. In the full cube $2^J$ we also have $e_j\to0$ along $\mathcal U$, since $\mathcal U$ is free. Composing $h$ with the inclusion into $2^J$, continuity and uniqueness of ultrafilter limits force $h(x)=0$, a contradiction.
\end{proof}

\begin{remark}
Proposition~\ref{prop:ultrafilter} gives an obstruction even when the number of source coordinates is arbitrary. For the equal-cardinality problem, classical compactness results yield a further obstruction (Corollary~\ref{cor:compact-obstruction} below); the absence of a free countably complete ultrafilter alone does not remove it. In Proposition~\ref{prop:countable-exceptions} the two uses of Lemma~\ref{lem:parameter-extension} depend on the fact that all proper stages of the construction are countable; simply running the same recursion through a larger cardinal does not establish the analogue.
\end{remark}

\begin{proposition}\label{prop:paramobstruction}
Let $P$ be a space, let $T$ be nonempty and separable, and let $Z\subset P$. If a continuous map $q:P\times T\to2$ satisfies $q(\{x\}\times T)=\{0\}$ for $x\in Z$ and $q(\{x\}\times T)=2$ for $x\notin Z$, then $Z$ is a $G_\delta$-subset of $P$.
\end{proposition}
\begin{proof}
For a countable dense set $\{t_n:n\in\N\}$ in $T$ (allow repetitions), continuity gives $Z=\bigcap_n\{x:q(x,t_n)=0\}$. The terms are clopen.
\end{proof}

In an uncountable product $\N^I$, no singleton is a $G_\delta$: the intersection of countably many basic neighborhoods of a point leaves some coordinate free. In particular, the assertion of Lemma~\ref{lem:parameter-extension} fails for general nonmetrizable zero-dimensional parameter spaces, even with the constant zero map prescribed on a singleton times $\B$. Composing a putative extension with the first binary coordinate would contradict Proposition~\ref{prop:paramobstruction}. This is a methodological obstruction: one cannot replace the Polish bases in the two extensions of Proposition~\ref{prop:countable-exceptions} by arbitrary nonmetrizable bases while retaining the same extension assertion. It does not rule out a different construction of a condensation at larger cardinals. The compactness obstruction below, in contrast, excludes the existence of a surjection itself for a specific target.

\begin{corollary}\label{cor:decomp}
If $\kappa$ is uncountable and $\varnothing\ne U\subseteq2^\kappa$ is open, then $U$ is a topological sum of a nonempty finite or countable family of open subspaces $U_n$, each $G_\delta$-dense in some clopen subspace homeomorphic to $2^\kappa$.
\end{corollary}
\begin{proof}
Choose $A$ by Lemma~\ref{lem:simultaneous-reflection}, applied to the constant sequence $U_n=U$, and put $V=\pi_A(U)$. Product projections are open, so $V$ is a nonempty open subset of the metrizable Cantor cube $2^A$. Enumerate a countable cover of $V$ by basic clopen cylinders contained in $V$; successively subtract the preceding finitely many cylinders and discard empty pieces. This partitions $V$ into nonempty pairwise disjoint sets $V_n$ clopen in $2^A$. Each is homeomorphic to the Cantor space. Put $K_n=V_n\times2^{\kappa\setminus A}$ and $U_n=U\cap K_n$, identifying the index cardinal with its underlying coordinate set. Each $K_n$ is clopen in $2^\kappa$ and homeomorphic to $2^\kappa$, since $V_n\cong\C=2^\N$ and $|\kappa\setminus A|=\kappa$.

To see that $U_n$ is $G_\delta$-dense in $K_n$, take a point $(y,z)$ in a nonempty $G_\delta$-subset $G$ of $K_n$. In each open set of a countable representation of $G$, choose a basic product neighborhood of $(y,z)$; the union of the finite sets of tail coordinates used in these neighborhoods is countable, say $J$. Since $y\in V_n\subseteq V$, the section $U_y$ is nonempty and thus $G_\delta$-dense. It contains a point $z'$ agreeing with $z$ on $J$. Then $(y,z')\in U_n\cap G$. Finally, the $U_n$ are pairwise disjoint, open in $U$, and cover it, so they give the asserted topological sum.
\end{proof}

\begin{remark}
The reflection lemmas do not require regularity of the index cardinal. For every uncountable $\kappa$, the preceding decomposition reduces the problem for all nonempty open subsets to the problem for $G_\delta$-dense open subsets of $2^\kappa$: combine condensations on the finitely or countably many summands, after partitioning one discrete coordinate of the source. For larger $\kappa$, $G_\delta$-density is weaker than negligibility of the complement in the sense of \cite{SV}, which excludes nonempty intersections of fewer than $\kappa$ open sets. The recursive part of the proof of Theorem~\ref{thm:main-gdelta}, however, uses the countability of every proper stage below $\omega_1$.
\end{remark}

\subsection{A compactness obstruction to surjectivity}
Recall that a space is finally $\kappa$-compact if every open cover has a subcover of cardinality strictly less than $\kappa$. The following result places Proposition~\ref{prop:extent} in the context of classical work on powers of discrete spaces. For related large-cardinal results for products of discrete factors, see Keisler and Tarski \cite[Theorem 4.32]{KT}. In particular, the results of Mr\'owka discussed by Lipparini \cite[p.~239]{Lipparini} imply that $\N^\kappa$ is finally $\kappa$-compact whenever $\kappa$ is weakly compact. Lipparini \cite[Theorem 5, pp.~242--243]{Lipparini} also characterizes final $\kappa$-compactness of $\N^\kappa$ for regular $\kappa$ by $(\kappa,\kappa)$-compactness of the infinitary language $\mathcal L_{\omega_1,\omega}$.

\begin{corollary}\label{cor:compact-obstruction}
Let $\kappa$ be an infinite cardinal. If $\N^\kappa$ is finally $\kappa$-compact, then there is no continuous surjection $\N^\kappa\to2^\kappa\setminus\{0\}$. In particular, this conclusion holds if $\kappa$ is weakly compact, or more generally if $\kappa$ is regular and $\mathcal L_{\omega_1,\omega}$ is $(\kappa,\kappa)$-compact.
\end{corollary}
\begin{proof}
If a closed discrete subset $D$ of cardinality $\kappa$ exists, choose for each $d\in D$ an open set $O_d$ meeting $D$ precisely at $d$. The cover consisting of $\N^\kappa\setminus D$ and all the $O_d$ has no subcover of cardinality less than $\kappa$. This contradicts final $\kappa$-compactness. Apply Proposition~\ref{prop:extent}. The two sufficient conditions for final $\kappa$-compactness are the classical implications just cited.
\end{proof}

Since the punctured cube is a nonempty open (hence $G_\delta$) subset, Corollary~\ref{cor:compact-obstruction} also answers the general problem of Section~1 negatively at every weakly compact cardinal. This conclusion is conditional on the existence of such a cardinal. If the least weakly compact cardinal exists, it is below every measurable cardinal, since every measurable cardinal is a limit of weakly compact cardinals; hence it admits no free countably complete ultrafilter (the least cardinal admitting one is measurable). This illustrates why the equal-cardinality compactness obstruction is not subsumed by Proposition~\ref{prop:ultrafilter}. Conversely, failure of final $\kappa$-compactness is only a necessary condition for surjectivity in this argument; we do not infer the closed-discrete-set condition from that failure in general.

\section{Open problems}
For the equal-cardinality questions, Proposition~\ref{prop:extent} gives a necessary condition for condensation onto the punctured cube, Corollary~\ref{cor:mycielski} verifies that condition below the first weakly inaccessible cardinal, and Corollary~\ref{cor:compact-obstruction} provides negative instances. The singleton-fibre refinement of Theorem~\ref{thm:singleton} does not imply bijectivity. The following questions are not resolved by those partial results.

\begin{problem}
For which uncountable cardinals $\kappa$ is $2^\kappa\setminus\{0\}$ a condensation of $\N^\kappa$? In particular, does the condition in Proposition~\ref{prop:extent} suffice for a condensation rather than just a surjection?
\end{problem}
\begin{problem}
For which cardinals $\kappa>\aleph_1$ satisfying the closed-discrete-set condition of Proposition~\ref{prop:extent} does $\N^\kappa$ condense onto every nonempty open subspace of $2^\kappa$? By Corollary~\ref{cor:decomp} it suffices to consider $G_\delta$-dense open subspaces.
\end{problem}

\begin{problem}
For which cardinals $\kappa>\aleph_1$ satisfying the closed-discrete-set condition of Proposition~\ref{prop:extent} does $\N^\kappa$ condense onto every nonempty $G_\delta$-subset of $2^\kappa$? In particular, does the conclusion of Theorem~\ref{thm:main-gdelta} hold when $\kappa=\aleph_2$?
\end{problem}

For general targets, a necessary condition for a condensation from $\N^\kappa$ is that every nonempty relatively open subset have cardinality $2^\kappa$: its inverse image is a nonempty open subset of $\N^\kappa$, which has that cardinality. Theorem~\ref{thm:main-gdelta} establishes sufficiency for the nonempty $G_\delta$ class, not necessity of that descriptive condition. The cardinality condition just noted is necessary but no sufficiency assertion is made here. A characterization of more general subspaces, including countable unions of $G_\delta$-subsets, and analogous questions for non-dyadic cubes remain outside the scope of this paper.

\clearpage
\appendix
\section{A verification guide to the transfinite construction}\label{app:verification}
This guide records the logical dependencies and the exact assertions to be checked in Proposition~\ref{prop:countable-exceptions}. It is a companion to the proof, not a replacement for it.

\subsection*{A.1. Inputs to the recursion}
\begin{enumerate}[label=\textup{(\arabic*)},leftmargin=*]
\item \emph{Reflection.} Lemma~\ref{lem:simultaneous-reflection} uses maximal families of pairwise disjoint finite remainders and a countable iteration to choose one countable set for all the open sets. Lemma~\ref{lem:enlargement-and-intersection} uses countably determined cylinders and the Baire theorem. Proposition~\ref{prop:polish-parameter-reduction} gives $V=V_0\times2^L$, with $V_0$ zero-dimensional Polish and $L$ countably infinite. Lemma~\ref{lem:parameter-absorption}, using Lemmas~\ref{lem:characterizations} and \ref{lem:basic-condensation}, gives $b:\B\to V$.
\item \emph{Normalization.} Lemma~\ref{lem:fresh-coordinate-witness} gives finite witnesses outside any specified countable set. In Proposition~\ref{prop:simultaneous-normalization}, the fresh blocks are chosen after collecting countably many witnesses, and the clopen partitions of $V$ may differ from block to block. Countability of each proper stage below $\omega_1$ is used here to ensure that the previously used coordinates form a countable set. The resulting homeomorphism gives closed sets satisfying \eqref{eq:zero-slices}.
\item \emph{Polish extension tools.} Lemma~\ref{lem:parameter-extension} uses Lemma~\ref{lem:finite-approximation} and nested countable clopen partitions of a metrizable base. Lemma~\ref{lem:exceptional-collapse} uses that $Z_\alpha$ is a nonempty closed subspace of a Polish space and that $Z_\alpha\times\B$ is homeomorphic to $\B$ for $\alpha<\omega_1$. Neither tool is applied with an uncountable-stage parameter.
\end{enumerate}

\subsection*{A.2. Preliminary identities and initial stage}
Fix the normalized sets $K_n$. The omitted factors of the prefix projections are compact. This gives closedness of $K_{n,\alpha}$ and, by compactness of the section at a \emph{fixed} parameter $v$, the limit identity \eqref{eq:K-limit}. No compactness of $V$ is needed. Definitions \eqref{eq:E-zero}--\eqref{eq:E-limit-definition} then give $K_{n,\alpha}\subseteq E_{n,\alpha}$, coherence \eqref{eq:E-coherence}, and equality \eqref{eq:E-equals-K-limit} at nonzero limit stages. At a successor stage the two sets need not be equal.

The domain factor $Z_\alpha$ in \eqref{eq:Z-definition} is closed and nonempty for every countable $\alpha$, including $0$, and its membership at limits is determined by all proper prefixes. It is important that nonemptiness is not the same as surjectivity of the restriction maps (Remark~\ref{rem:nonextendible}). At $\alpha=0$, $H_0=b$ under the natural identifications, $Z_0$ is a singleton, and $E_{n,0}=V$. Thus \eqref{eq:fibre-product}--\eqref{eq:L-definition} hold even if some $K_{n,0}$ is empty.

\subsection*{A.3. Successor step: check the whole fibre}
Fix $\alpha<\omega_1$ and abbreviate $E=E_{n,\alpha}$, $K=K_{n,\alpha}$ for one index $n$. Both are closed in the zero-dimensional Polish space $Y_\alpha$. The two successive applications of Lemma~\ref{lem:parameter-extension} have the following data:
\begin{center}
\begin{tabular}{@{}lll@{}}
\textit{Base} & \textit{Closed prescribed subset} & \textit{Variable of condensation} \\
$E$ & $K$ & $(x,u)\in Z_\alpha\times\B$ \\
$Y_\alpha\times P_\alpha$ & $E\times Z_\alpha$ & $u\in\B$
\end{tabular}
\end{center}
The first extends $r_\alpha$ and the second extends the resulting $R_{n,\alpha}$. The latter agrees with the former on its entire prescribed domain, not just on $K$. All bases and prescribed subsets satisfy the hypotheses; the empty cases are allowed.

For a desired output $y^+$ let $y=\pi_\alpha^{\alpha+1}(y^+)$ and $d=z_n^+(\alpha)$. The primary inputs are uniquely determined by $b$ and $c$. At fixed $y$, the $n$th new equation constrains only $(x_n,u_n)$. The possibilities are:
\begin{center}
\begin{tabular}{@{}lll@{}}
\textit{Old output and new value} & \textit{Allowed pairs} & \textit{New exceptional status} \\
$y\notin E$ & unique & outside $E_{n,\alpha+1}$ \\
$y\in E\setminus K$ & unique & outside $E_{n,\alpha+1}$ \\
$y\in K,\ d=0$ & exactly $Z_{\alpha+1}$ & inside $E_{n,\alpha+1}$ \\
$y\in K,\ d\ne0$ & unique & outside $E_{n,\alpha+1}$
\end{tabular}
\end{center}
In the first row, the old $x_n$ is fixed and the second extension is bijective in $u_n$. In the second row, the entire pair is determined by the first extension. In the last two rows use the zero and nonzero fibres of $r_\alpha$, respectively.

For both inclusions of the new fibre identity, necessity follows from prefix compatibility and the new equations. Sufficiency uses the full old fibre identity: any family of allowed pairs has old prefix in $H_\alpha^{-1}(y)$, hence its old output really is $y$, the parameter used in every new equation. For example, the diagonal of $\{0,1\}^2$ projects onto both factors but does not contain all combinations; projection identities alone would not justify this step. Continuity is checked coordinatewise in \eqref{eq:H-successor}, independently of the fibre calculation. All factors are nonempty, so $H_{\alpha+1}$ is onto.

\subsection*{A.4. Limit stages and termination}
At a nonzero limit $\lambda$, each output block is already specified by an earlier compatible map. This defines a continuous $H_\lambda$ on the full product, whether or not its image is all of $Y_\lambda$. At a countable limit, substitute \eqref{eq:fibre-product} into \eqref{eq:fibre-limit}. The auxiliary constraints are
\[
 \forall\alpha<\lambda\ \forall n\in\N:\quad
 x_n\res\alpha\in L_{n,\alpha}(y_\alpha).
\]
Interchanging the two universal quantifiers gives exactly the factors in \eqref{eq:M-factor}. If the $n$th exception persists, that factor is $Z_\lambda$ by \eqref{eq:Z-limit}. If it does not persist, coherence of $E_{n,\alpha}$ guarantees that all later factors are singletons. For any two such stages, choose a point in the later full fibre (which is nonempty by the induction hypothesis) and restrict it to the earlier stage. This proves compatibility of the singleton tuples, including the required earlier restrictions. Their union is the unique admissible tuple. At no stage is surjectivity of the maps $Z_\delta\to Z_\alpha$ needed.

At $\lambda=\omega_1$ the same compatibility defines the continuous map. If $y\in K_n$, all its prefixes lie in $E_{n,\alpha}$ and a preimage would yield an element of $Z_{\omega_1}$, contradicting \eqref{eq:Z-empty}. If $y\notin\bigcup_nK_n$, equality \eqref{eq:E-equals-K-limit} yields a countable exit stage $\beta_n$ for each $n$. Since $\operatorname{cf}(\omega_1)>\omega$, their supremum $\beta$ is countable. All full fibres at countable stages at or above $\beta$ are singletons. The union of their compatible tuples is a preimage of $y$, and every preimage must have the same cofinal prefixes. This proves \eqref{eq:final-fibres}. Finally the domain has $\aleph_1$ countable discrete factors and is homeomorphic to $\N^{\omega_1}$.

The obstruction in Proposition~\ref{prop:paramobstruction} concerns the extension tool used in the successor step; it is not a counterexample to a condensation at $\aleph_2$. By contrast, Corollary~\ref{cor:compact-obstruction} excludes even a surjection onto a particular open target at weakly compact cardinals. These limitations have different logical roles.
\Needspace{8\baselineskip}
\section*{AI assistance} ChatGPT (OpenAI) was used to explore research questions, develop mathematical constructions and proof arguments, identify relevant literature, draft and revise the text, and assist in checking proofs and references. These AI-assisted checks do not constitute independent expert verification. The author takes responsibility for the final mathematical content and bibliographic accuracy.

\end{document}